\documentclass[journal]{IEEEtran}

\usepackage{cite}
\usepackage{todonotes}
\usepackage{graphicx}
\usepackage{subfigure}
\usepackage{epstopdf}
\usepackage{url}
\usepackage{stfloats}
\usepackage{latexsym}
\usepackage{amsfonts}
\usepackage{amsmath}
\usepackage{amssymb}
\usepackage{array}
\usepackage{textcomp}
\usepackage{pxfonts}
\usepackage{txfonts}
\usepackage{cases}
\usepackage{psfrag}
\usepackage{epsfig}
\usepackage{setspace}
\usepackage{enumitem}
\usepackage{xcolor}
\usepackage{multirow}
\usepackage{multicol}
\usepackage{lmodern}

\usepackage[T1]{fontenc}
\usepackage{color}
\usepackage{float}
\usepackage{esint}
\usepackage[unicode=true,pdfusetitle, bookmarks=true,bookmarksnumbered=false,bookmarksopen=false,  breaklinks=false,pdfborder={0 0 0},backref=false,colorlinks=true]{hyperref}

\makeatletter

\floatstyle{ruled}
\newfloat{algorithm}{tbp}{loa}
\providecommand{\algorithmname}{Algorithm}
\floatname{algorithm}{\protect\algorithmname}

\usepackage{fullpage}
\usepackage{algorithmic,algorithm}

\makeatother

\newtheorem{Definition}{Definition}
\newtheorem{Remark}{Remark}

\usepackage{todonotes}   

\newcommand{\R}{\ensuremath{\mathbb{R}}}
\newcommand{\N}{\ensuremath{\mathbb{N}}}

\newlist{aims}{enumerate}{1}
\setlist[aims,1]{
  label={Additive Faults, set~\arabic*:},
  leftmargin=*,
  align=left,
  labelsep=10mm,
}

\newlist{steps}{enumerate}{1}
\setlist[steps,1]{
  label={Step~\arabic*:},
  leftmargin=*,
  align=left,
  labelsep=10mm,
}

\begin{document}

\title{\Large Bayesian inference in non-Markovian state-space models with applications to fractional order systems}
\author{Pierre E. Jacob, S.M.Mahdi Alavi, Adam Mahdi, Stephen J. Payne and David A. Howey
\thanks{P.E. Jacob is with the Department of Statistics, Harvard University, Science Center 7th floor, 1 Oxford Street, Cambridge, MA 02138-2901, USA.  Email: pierre.jacob.work@gmail.com}
\thanks{S.M.M. Alavi was with the Energy and Power Group, Department of Engineering Science, University of Oxford. He is now with the Brain Stimulation Engineering Laboratory, Department of Psychiatry \& Behavioral Sciences, Duke University, Durham, NC 27710, USA. Email: mahdi.alavi@duke.edu}
\thanks{A. Mahdi and S.J. Payne are with the Institute of Biomedical Engineering, Department of Engineering Science, University of Oxford, Old Road Campus Research Building, Oxford, OX3 7DQ, United Kingdom. Emails:  \{adam.mahdi, stephen.payne\}@eng.ox.ac.uk}
\thanks{D.A. Howey is with the Energy and Power Group, Department of Engineering Science, University of Oxford, Parks Road, Oxford, OX1 3PJ, United Kingdom. Email: david.howey@eng.ox.ac.uk}
\thanks{The source code is available online at \url{https://github.com/pierrejacob/BatteryMCMC}}
}

\markboth{ }%
{Shell \MakeLowercase{\textit{et al.}}: Bare Demo of IEEEtran.cls
for Journals}
\maketitle

\begin{abstract}
Battery impedance spectroscopy models are given by fractional order (FO) differential
equations. In the discrete-time domain, they give rise to state-space models
where the latent process is not Markovian. Parameter estimation for these models
is therefore challenging, especially for non-commensurate FO models.
In this paper, we propose a Bayesian approach to identify the parameters
of generic FO systems. The computational challenge is tackled with particle 
Markov chain Monte Carlo methods, with an implementation specifically designed
for the non-Markovian setting. The approach is then applied to estimate
the parameters of a battery non-commensurate FO equivalent circuit model.
Extensive simulations are provided to study the practical identifiability of
model parameters and their sensitivity to the choice of prior distributions,
the number of observations, the magnitude of the input signal and the
measurement noise. 
\end{abstract}


\IEEEpeerreviewmaketitle

\section{Introduction}
\label{sec:introduction}
Fractional Order (FO) models are important in the study of electrochemical and biological systems, \cite{Monje2010} and references therein. Certain features in the FO models make their identification challenging. Before reviewing specific issues, we recall that a continuous-time state-space model of FO systems is given by \cite{Alavi2015b}:
\begin{eqnarray}\label{SS:c}
\begin{aligned}
\frac{d^{\alpha} x(t)}{dt^{\alpha}}&= \bar{A}(\beta)\, x(t) + \bar{B}(\beta)\,u(t),\\
 y(t) &= M(\beta)\, x(t)+ D (\beta)\,u(t),
\end{aligned}
\end{eqnarray}
where $x \in \R^n$ is the state vector; $u \in \R$ and $y \in \R$ are input and output signals, respectively; $\bar{A}(\beta)\in \R^{n \times n}$, $\bar{B}(\beta) \in \R^{n \times 1}$, $M(\beta) \in \R^{1 \times n}$ and $D(\beta) \in \R$ are system matrices which depend on the parameter vector $\beta \in \R^q$ to be identified. Moreover,
\begin{equation}\label{fracvec}
\frac{d^{\alpha} x(t)}{dt^{\alpha}} = \Big[\frac{d^{\alpha_1} x_1(t)}{dt^{\alpha_1}},\ldots,\frac{d^{\alpha_n} x_n(t)}{dt^{\alpha_n}}\Big]
\end{equation}
is the vector of FO derivatives, with unknown FOs  denoted by $\alpha_i\in(0,1)$, $i=1,\ldots,n$. 

By defining the parameter vector as
\begin{align}
  \label{parvec}
  \theta=\left[\begin{array}{cccc} \alpha_1 & \cdots & \alpha_n &  \beta \end{array}\right],
\end{align}
the corresponding model in the discrete-time domain is given by \cite{Alavi2015b}:  
\begin{eqnarray}\label{SS:d}
\begin{aligned}
 x_{k+1}&=\displaystyle \sum_{j=0}^{k} A_j(\theta)\, x_{k-j}+B(\theta) u_k, \\
\label{output equation discrete}
y_k&=M(\theta) x_k + D(\theta) u_k,
\end{aligned}
\end{eqnarray}
with
\begin{eqnarray}\label{SS:d-matrices}
\begin{aligned}
& A_0(\theta)=\mbox{diag}\{\alpha_1,\cdots,\alpha_n\}+\\
& \hspace{3em}\mbox{diag}\{T_s^{\alpha_1},\cdots,T_s^{\alpha_n}\}\bar{A}(\beta),\\
& A_j(\theta)=(-1)^{j}\mbox{diag}\left\{\binom{\alpha_1}{j+1},\cdots,\binom{\alpha_n}{j+1}\right\},\\
&\hspace{3em} \mbox{~for~} 1 \leq j\\
& B(\theta) =\mbox{diag}\{T_s^{\alpha_1},\cdots,T_s^{\alpha_n}\}\bar{B}(\beta),
\end{aligned}
\end{eqnarray}
where $T_s$ is the sample time, $k\in \N$ is the time index,  $\mbox{diag}\{\cdot\}$ denotes the diagonal matrix and $\binom{\alpha_i}{j}$ is the binomial coefficient given by
\begin{equation}\label{Gamma_func}
\binom{\alpha_i}{j}=\frac{\Gamma(\alpha_i+1)}{\Gamma(j+1)\Gamma(\alpha_i+1-j)},
\end{equation}
where, $\Gamma(\cdot)$ denotes the gamma function
\[
\Gamma(\alpha_i)=\int_{0}^{\infty} z^{\alpha_i-1} e^{-z}dz,~ \text{for} ~ \alpha_i\in\mathbb{C}~ \text{with}~  \Re(\alpha_i)>0,
\]
where $\Re$ denotes the real part of a complex number.

\begin{Definition}
\label{Def:commensurateFOsys}An FO system is said to be commensurate if for all
$i\in\{1,\ldots n\}$, there exists $\rho \in \N,~\mbox{such~that}~\alpha_i =
\rho \alpha$, where $\alpha \in \R$; otherwise it is said to be
non-commensurate \cite{Monje2010}.
\hfill{$\square$}
\end{Definition}


Now, the main issues associated with parameter estimation in FO
systems are more evident: 1) the state vector $x_{k+1}$ depends on all the past
states $x_0$ up to $x_k$, therefore the FO system is non-Markovian \footnote{In
Markovian systems, $x_{k+1}$ can be written as functions of $x_k$ and inputs.},
and 2) the model is highly nonlinear, with respect to the parameters. 

Several least-squares estimation methods have been proposed in \cite{LeLay1998, Lin2001, Cois2002} for the identification of continuous-time FO transfer functions of the form
\begin{align}
\label{General-Fractional-order-sys-TF}
F(s)=\frac{Y(s)}{U(s)}=\frac{b_m s^{\beta_m}+\cdots+b_1s^{\beta_1}+b_0s^{\beta_0}}{a_n s^{\alpha_n}+\cdots+a_1s^{\alpha_1}+1},
\end{align}
where $U$ and $Y$ are Laplace transforms of the input and output (observation) signals, respectively.
These methods have been modified in \cite{Malti2008, Victor2013}, and developed into the Crone toolbox \cite{Oustaloup2000, Melchior2002}. The 
Crone toolbox is mainly based on the instrumental variable state variable filter (ivsvf) method and the simplified refined instrumental variable for continuous-time fractional (srivcf) method, which are both based on instrumental variable concepts \cite{Young1979, Jakeman1979, Young1980a}. In these methods, the FO model is linearised by approximating all fractional differentiation operators $s^{\alpha_i}$ and $s^{\beta_j}$ by higher-order rational transfer functions \cite{Malti2008}. The coefficients $b_j$'s and $a_i$'s are identified by using the coefficient map that exists between the original and approximated models. Manual search \cite{Malti2008}, gradient descent \cite{Victor2013}, or interior-point \cite{Alavi2015} optimizations are combined with the ivsvf and srivcf functions for the estimation of the fractional orders. In the non-commensurate case, the approximate model is of high order so that the coefficient map between the original and approximated model becomes very complex and may be intractable. This issue highlights the need for novel tools to directly identify general FO models.

In \cite{Wang2015}, an identification method based on swarm optimization has been proposed to identify a battery non-commensurate FO model.  

In this paper, a novel method based on Bayesian inference is presented.
Bayesian inference assumes a prior distribution on parameters, which is then updated, using the observations, to yield the posterior distribution. By investigating the posterior distribution, parameters can be estimated, along with associated credible regions reflecting uncertainty. By comparing the prior and posterior distributions, the identifiability of the model can be assessed from a practical perspective, by answering questions such as: do the data inform the parameters? Are some parameters more easily identifiable than others?

Markov chain Monte Carlo methods provide an approach to approximate the posterior distribution,
provided that the associated probability density function can be evaluated point-wise. More recently, particle Markov chain Monte Carlo (PMCMC) \cite{Andrieu2010} methods have been proposed for Bayesian inference in state-space models, where the posterior density can only be estimated point-wise. For FO models, the added difficulty lies in the fact that the latent process is non-Markovian. Here we will use PMCMC, with a novel implementation that leads to practical approximations in the case of non-Markovian models. The method enables Bayesian inference without requiring any model simplifications, such as linearizations or Gaussian assumption, and is applicable to any state-space model with non-Markovian latent processes.

\begin{figure}
\centering
\includegraphics[scale=.8]{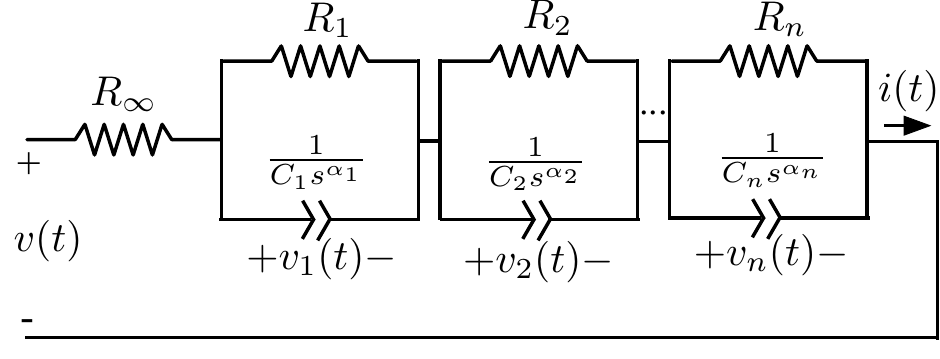}
\caption{The general battery electrochemical impedance spectroscopy model.}
\label{Fig:EISgeneral}
\end{figure}

Recently, there has been a significant interest in the design of model-based battery systems to improve the efficiency and reliability of electric vehicles and renewable energies \cite{Howey2015}. Among the employed models, the battery Electrochemical Impedance Spectroscopy (EIS) models with FO dynamics have received much attention. They are more accurate than the conventional lumped models and are computationally less expensive than the electrochemical models defined by highly coupled partial differential equations. A comprehensive survey of battery models has been given in \cite{Howey2015}. A general schematic of the battery impedance FO models is shown in Figure \ref{Fig:EISgeneral}, where $i$ and $v$ denote the battery current and voltage, respectiely. $R_{\infty}$ represents the battery ohmic resistance at high frequencies. Each parallel pair is employed to model the battery processes over a certain frequency range. The number of parallel pairs depends on the required accuracy for the frequency domain fitting of impedance spectra. The terms $C_is^{\alpha_i}$, $i=1,\cdots,n$, called constant phase elements (CPEs), model diffusion processes (e.g. charge transfer resistance and double layer capacitance) more accurately compared to the lumped models as shown in \cite{Barsoukov2005} and \cite{Alavi2015}. In low frequency ranges, the impedance frequency response may show constant phase behaviour such that the associated parallel resistor can be considered as an open circuit. This is referred to as Warburg term in the literature \cite{Barsoukov2005}.  Reference \cite{Barsoukov2005} provides more information on battery EIS and associated FO models. 

By defining $u=i$ and $y=v$, and the voltage across the CPEs as the state variables,
\begin{align}
x \triangleq \left[\begin{array}{ccc} v_1 & \cdots & v_n\end{array}\right],
\end{align}
it is easy to show that $A_j$, $B$, $M$ and $D$ in the state-space model \eqref{SS:d} are given by:
\begin{eqnarray}\label{gen-ecm-ssd}
\begin{aligned}
& A_j(\theta)=\mbox{diag}\left\{a_{1,j}(\theta), \cdots, a_{n,j}(\theta)\right\}\\
& B(\theta)=\left[\begin{array}{ccc} b_1(\theta) & \cdots & b_n(\theta)\end{array}\right]\\
& M(\theta)=\left[\begin{array}{ccc} m_1 & \cdots & m_n\end{array}\right]\\
& D(\theta)=d(\theta),
\end{aligned}
\end{eqnarray}
with
\begin{eqnarray}\label{gen-ecm-params}
\begin{aligned}
& a_{i,0}(\theta)= \alpha_i-\frac{T_s^{\alpha_i}}{R_iC_i},~ a_{i,j}(\theta)= (-1)^{j}\binom{\alpha_i}{j+1}\\
& b_i(\theta)=\frac{T_s^{\alpha_i}}{C_i},~ m_i=1,~d(\theta)=R_\infty\\
& \text{for~}i=1,\ldots,n \mbox{~and~} j=1,2,\ldots,T,
\end{aligned}
\end{eqnarray}
where $T$ is the data length, $k=1,\ldots,T$. 
 
In order to estimate the parameter vector $\theta$ from $v$ and $i$, the battery FO models are currently fitted to frequency domain impedance spectra that are obtained through EIS, \cite{Troeltzsch2006, Macdonald1982, Boukamp1986}. This requires a data conversion that may introduce bias in the estimation \cite{Alavi2015}. Thus, parameter estimation of the battery FO models directly from time-domain data is very appealing.

In \cite{Alavi2015b} we showed that the FO model with just a single CPE (where $n=1$ in Figure \ref{Fig:EISgeneral}) is globally identifiable. Here we study the identifiability of the model with more than one CPE using the proposed Bayesian approach. Extensive simulations on a battery FO system are provided to illustrate how the proposed method enables the study of various effects on parameter identification, such as the data length, the magnitude of the input signal, the choice of prior, and the measurement noise.  We also show that the parameters of the Warburg term in the non-commensurate FO battery model are not identifiable. The source code is available online at
\url{https://github.com/pierrejacob/BatteryMCMC}.

The main contributions of the paper are:
\begin{itemize}
  \item[-] Introduction of Bayesian framework for Identification of general FO systems.
  \item[-] Reduction of the computational and memory cost for the non-Markovian particle filter.
  \item[-] Applications to identifiability of non-commensurate FO battery models. 
\end{itemize}


\section{Bayesian Inference}
\label{sec:pmcmc}
In this section, we describe Bayesian inference in the context of non-Markovian state-space models.
The vector $(v_{s},\ldots,v_{t})$ (resp. $(v^{s},\ldots,v^{t})$) is denoted by
$v_{s:t}$ (resp. $v^{s:t}$). The notation $a\sim p(\cdot\mid b)$ means that $a$
is a random variable distributed according to a distribution $p$ which depends
on $b$. 
The Gaussian distribution with mean $\mu$ and variance $\sigma^{2}$ is denoted by $\mathcal{N}(\mu,\sigma^{2})$.

\subsection{Inference in state-space models\label{sub:statespacemodels}}
Given measurements $y_{0:T}$, parameter inference refers to the task of estimating the parameter vector $\theta$ in
the general state-space model
\begin{eqnarray}\label{bayes-model}
\begin{aligned}
& x_{0}\sim\mu(\cdot\mid\theta),\\
    \forall k\in \{1,\ldots,T\}\;\;\;& x_{k}\sim f_{k}(\cdot\mid x_{0:k-1},\theta),\\
\forall k\in \{0,\ldots,T\}\;\;\;& y_{k}\sim g_{k}(\cdot\mid x_{k},\theta),
\end{aligned}
\end{eqnarray}
where  $T\in\mathbb{N}$ denotes the total number of observations
\cite{Cappe2005, Wasserman2011}.  Note that the notation reflects that
$f_{k}$ potentially depends on all the past states $x_{0:k-1}$, and thus
$\{x_0, \cdots, x_{k}\}$ is not necessarily a Markov chain. Both $f_{k}$ and $g_{k}$ could be nonlinear functions. Therefore, the FO model \eqref{SS:d} is in the Bayesian framework
\eqref{bayes-model}, upon specifying a state noise and an observation noise. For simplicity,
we will consider additive Gaussian noises, that is, we consider:
\begin{eqnarray}\label{SS:d:withnoise}
\begin{aligned}
 x_{k+1}&=\displaystyle \sum_{j=0}^{k} A_j(\theta)\, x_{k-j}+B(\theta) u_k +\sigma_x \varepsilon_k, \\
y_k&=M(\theta) x_k + D(\theta) u_k + \sigma_y \eta_k,
\end{aligned}
\end{eqnarray}
where $(\varepsilon_k)$ and $(\eta_k)$ are sequences of independent standard Gaussian variables, 
and $\sigma_x$, $\sigma_y$ are positive values.
Statistical inference classically relies on the likelihood function, defined by $\theta\mapsto p(y_{0:T}\mid\theta)$, where $p(y_{0:T}\mid\theta)$ is
the density of the observations evaluated with the dataset $y_{0:T}$ for the parameter $\theta$. For state-space models the likelihood can be written, for any $\theta$, in terms of $\mu$, $f_{k}$ and $g_{k}$ as follows,
\begin{align}
\nonumber
 p(y_{0:T}\mid\theta)&= \int \mu\left(x_{0}\mid\theta\right) g_0(y_0 \mid x_0,\theta) \times \\ & \prod_{k=1}^{T}f_{k}\left(x_{k}\mid x_{0:k-1},\theta\right)g_{k}(y_{k}\mid x_{k},\theta)dx_{0:T}.\label{eq:likelihood}
\end{align}
One way to estimate the parameters is to maximise the likelihood function
$\theta\mapsto p(y_{0:T}\mid\theta)$ with respect to $\theta$, yielding the
maximum likelihood estimator. Confidence intervals can be constructed by
assuming that the maximum likelihood estimator is asymptotically normal;
other ways to obtain confidence intervals, such as bootstrap, are not readily
applicable to state-space model settings.
Alternatively, Bayesian inference relies on
a probability distribution defined on the space $\Theta$, representing our
knowledge about the parameters given the observations \cite{Robert2007}. We
first specify a prior distribution $p(\theta)$, representing the knowledge of
$\theta$ before observing the data (for instance, based on past experiments,
intuition, literature search, etc). Then, the posterior distribution of the
parameters given the observations is given by Bayes formula:
\begin{equation}
p(\theta\mid y_{0:T})=\frac{p(\theta)p(y_{0:T}\mid\theta)}{\int p(\vartheta)p(y_{0:T}\mid\vartheta)d\vartheta}.\label{eq:bayesformula}
\end{equation}
In the numerical experiments, we will be particularly interested in the
changes between the prior and the posterior, which inform us about how much
has been learned from the data, and thus give a practical notion of identifiability.
The posterior distribution $p(\theta\mid
y_{0:T})$ can also be used to obtain point estimates, such as the posterior
mean, and uncertainty is typically measured using credible regions
\cite{Robert2007}. Since the posterior distribution typically does not belong
to a standard class of parametric probability distributions, it is approximated
using Monte Carlo methods, which provide samples approximately distributed 
according to the posterior distribution.

For general state-space models, that is, generic choices of $\mu$, $f_{k}$ and
$g_{k}$, the integral in Eq. (\ref{eq:likelihood}) cannot be evaluated exactly,
and has to be numerically approximated. Therefore, exact posterior density
evaluations are not available either. We will describe how to construct
flexible and efficient approximations of the posterior distribution using a
combination of Markov chain Monte Carlo (MCMC) algorithms \cite{Robert2013} and
particle filters \cite{Doucet2001}, called particle Markov chain Monte Carlo
(PMCMC) \cite{Andrieu2010}.

\subsection{Likelihood estimation using particle filters \label{sub:particlefilters}}
For our purposes, particle filters \cite{Gordon1993,
Doucet2001, Doucet2009} will provide
approximations of the likelihood in Eq. \eqref{eq:likelihood}, 
for any $\theta\in\Theta$, up to a multiplicative constant.
The user sets a number of particles $N\in\mathbb{N}$. Larger values of $N$ yield better
precision but the computational cost of the method increases linearly with $N$.
A population of $N$ particles is first sampled from the initial
distribution $\mu(\cdot\mid\theta)$. These particles are denoted by
$x_{0}^{1:N}=(x_{0}^{1},\ldots,x_{0}^{N})$. Then, each particle is weighted,
using the first observation $y_{0}$. That is, for each $i\in\left\{ 1,\ldots,N\right\}
$, $w_{0}^{i}=g_{0}(y_{0}\mid x_{0}^{i},\theta)$. These weights $w_{0}^{1:N}$
are such that the weighted sample, $(w_{0}^{i},x_{0}^{i})_{i=1}^{N}$ is
approximately distributed according to $p(x_{0}\mid y_{0},\theta)$. 
This weighting is commonly called importance
sampling \cite{Robert2013} and concludes the initialization of the algorithm.

The next step consists in selecting some particles and discarding others,
according to their weights. This is the resampling step, which comes in various flavors.
The most standard resampling scheme, called multinomial resampling,
consists in drawing, independently for each $i\in\left\{ 1,\ldots,N\right\} $,
an ancestor variable $a_{0}^{i}\in\left\{ 1,\ldots,N\right\} $ according to
a categorical distribution with parameters
$(w_{0}^{1}/\sum_{j=1}^{N}w_{0}^{j},\ldots,w_{0}^{N}/\sum_{j=1}^{N}w_{0}^{j})$.
The variable $a_{0}^{i}$ is interpreted as the index of the parent of particle
$i$ at time $1$, among the $N$ particles at time $0$. 
Once the parents $a_{0}^{1:N}$ are chosen, the new particles are sampled and weighted according to
$ x_{1}^{i}\sim q_{1}(\cdot\mid x_{0}^{a_{0}^{i}},\theta)$
and
\begin{align*}
w_{1}^{i}=\frac{f_{1}(x_{1}^{i}\mid x_{0}^{a_{0}^{i}},\theta)g_{1}(y_{1}\mid x_{1}^{i},\theta)}{q_{1}(x_{1}^{i}\mid x_{0}^{a_{0}^{i}},\theta)}.
\end{align*}
This is another importance sampling step, where the proposal distribution is
$q_{1}$, and the weights are computed such that
$(x_{1}^{i},w_{1}^{i})_{i=1}^{N}$ is approximately distributed according to
$p(x_{1}\mid y_{0},y_{1},\theta)$. Let us introduce the paths
$\bar{x}_{0:1}^{i}=(x_{0}^{a_{0}^{i}},x_{1}^{i})$, for each $i\in\left\{
1,\ldots,N\right\} $. Then the weighted paths
$(w_{1}^{i},\bar{x}_{0:1}^{i})_{i=1}^{N}$ are approximately distributed
according to the path distribution $p(x_{0},x_{1}\mid y_{0},y_{1},\theta)$. The
algorithm then proceeds in a similar fashion for the subsequent steps,
resampling, sampling and weighting the particles until all the data have been
assimilated. A pseudo-code description of the algorithm is given in Section
\ref{sec:algo:particlefilter}. The algorithm provides at every time $k$ a
weighted sample $(w_{k}^{i},\bar{x}_{0:k}^{i})_{i=1}^{N}$ that approximates the
path distribution $p(x_{0:k}\mid y_{0:k},\theta)$. More importantly for us, the quantity
\begin{equation}
\hat{p}(y_{0:T}\mid\theta)=\prod_{k=0}^{T}\left(\frac{1}{N}\sum_{i=1}^{N}w_{k}^{i}\right)\label{eq:likelihoodestimator}
\end{equation}
is an estimator of the likelihood $p(y_{0:T}\mid\theta)$. A rich literature is
devoted to the theoretical study of the algorithm and the properties of this
type of estimator \cite{DelMoral2004,Cerou2011,Whiteley2013}, at least in the settings
of Markovian latent processes. These results
indicate that, as the number of observations $T$ goes to infinity, the relative
variance of the likelihood estimator $\hat{p}(y_{0:T}\mid\theta)$ can be
bounded independently of $T$ if one chooses $N$ proportionally to $T$.

The user has to choose the proposal distributions $q_{k}(\cdot\mid
x_{0:k-1},\theta)$ for all $k\geq1$. The default choice of proposal is to use
$q_{k}=f_{k}$, the model transition. This choice leads to a simplified
expression of the weight, $w_{k}^{i}=g_{k}(y_{k}\mid x_{k}^{i},\theta)$.
Another choice, called locally optimal proposal, consists in using 
\begin{equation}
q_{k}(x_{k}\mid x_{0:k-1},\theta)=p(x_{k}\mid y_{k},x_{0:k-1},\theta),\label{eq:locallyoptimalproposal}
\end{equation}
which leads to the weight $w_{k}^{i}=p(y_{k}\mid\bar{x}_{0:k-1}^{i},\theta)$.
The locally optimal proposal takes the next observation $y_{k}$ into
account when propagating the particles from time $k-1$ to time $k$. 
Appendix \ref{sec:localproposal} details the calculation leading to 
the locally optimal proposal for the battery models of Section \ref{sec:batpmcmc}.
The choice of proposal distributions has a direct impact on the variance of the
estimator $\hat{p}(y_{0:T}\mid\theta)$. 

%
\subsection{Special features of non-Markovian models \label{sec:nonmarkovian}}

We consider the following implementation of the standard particle filter. Given the non-Markovianity of the latent
process, sampling each particle $x_{k}^{i}$ at time $k$ requires
computing a function $\psi_{k}$ of a trajectory, $(x_{0}^{i},\ldots,x_{k-1}^{i})$.
In the models considered in the article, such as \eqref{SS:d:withnoise}, this function takes the form
of a weighted sum:
$\psi_{k}(\bar{x}_{0:k}^{i})=\sum_{t=0}^{k}\alpha_{k,t}\bar{x}_{k-t}^{i}$,
where $(\alpha_{k,t})_{t\geq0}$, for each $k$, are coefficients that 
can be computed given the parameter $\theta$.
Naively computing this weighted sum for each particle would yield
a cost of order $k$, at time $k$; and thus an overall computational cost of order $N\times T^{2}$,
where $T$ is the number of observations to assimilate and $N$ the
number of particles. Furthermore, the naive memory cost of storing all the trajectories would be
of order $N\times T$.

However, we can reduce both the memory cost and the computational cost, 
by representing the trajectories as branches of a tree. 
At each
step of the particle filter, new leaves are added to the tree, corresponding
to the new generation of particles. The ancestors give the list of new branches.
At the same time, some existing branches
of the tree can be cut, corresponding to  particles that have been
discarded in the resampling step. As a result, it was shown in \cite{Jacob2013}
that the number of nodes in the tree is of order $k+CN\log N$ at time $k$, in expectation.
The constant $C$ does not depend on $N$ and $T$.
Efficient algorithms to implement this tree structure are given in \cite{Jacob2013}.
Therefore, the memory cost can be reduced from $N\times T$ to $T + C N\log N$.
Furthermore, computing the $N$ sums 
at time $k$ only requires browsing the whole tree once, from the root to the leaves, for a cost
of order $k+CN\log N$ at step $k$. Performing this operation for each $k$ in
$\left\{ 1,\ldots,T\right\} $ yields an overall computational cost of order $T^{2}+CTN\log N$,
instead of $N\times T^2$ with the standard implementation..

Thus, the memory cost of particle filtering for the paths of a non-Markovian
model is of order $T+CN\log N$, and the computational cost of the
algorithm is of order $T^{2}+CTN\log N$, when the transition involves
sums as in Eq. \eqref{SS:d:withnoise}. The same would be true 
for any calculation that only requires parsing the tree structure once.
This is to be compared with a computational cost of $NT$ and a memory cost of $N$ for Markovian models. 
Since $N$ is typically comparable to $T$, and $\log N$ is a small value, the tree structure gives an order of magnitude improvement, both computationally and in terms of memory, over the standard implementation.

\subsection{Parameter estimation using particle Markov chain Monte Carlo \label{sub:pmcmc}}

We next describe briefly how point-wise likelihood estimates
$\hat{p}(y_{0:T}\mid\theta)$, such as the ones produced by particle filters,
can be used to obtain samples from the posterior distribution. 
We consider MCMC algorithms to produce $\theta_{1:M}=(\theta_{1},\ldots,\theta_{M})$,
the realization of a Markov chain approximately following the posterior distribution 
of Eq. \eqref{eq:bayesformula} \cite{Robert2013}. 
Standard MCMC requires point-wise evaluations of the posterior density,
and thus is not applicable here. Instead, we have access to estimators provided by particle filters,
which leads to particle MCMC \cite{Andrieu2010}.
These methods have been thoroughly analyzed in
\cite{Andrieu2015,Andrieu2013,Lindsten2014,Chopin2015}. We will use in particular 
the Particle Marginal Metropolis--Hastings (PMMH) algorithm, based
on the standard Metropolis--Hastings algorithm \cite{Liu2008,Robert2013}. 
Remarkably, the Markov
chain $\theta_{1:M}$ produced by PMMH still converges to the target distribution when $M$ goes
to infinity, as if exact density evaluations were used. The PMMH algorithm is
described in pseudo-code in Section \ref{sec:algo:pmmh}, along with details on
algorithmic tuning.

\begin{Remark}
The literature on particle filters and PMCMC methods usually considers the case
where the latent process $(x_{k})_{k\in\mathbb{N}}$ is a Markov chain, but we see that
all the above algorithms are directly implementable in the non-Markovian setting.
Few articles have considered the non-Markovian case \cite{Chopin2011,Lindsten2012,Wood2014}. \hfill{$\square$}
\end{Remark}


\section{Bayesian Inference in Battery Systems}
\label{sec:batpmcmc}

\begin{figure}
\centering
\includegraphics[scale=0.8]{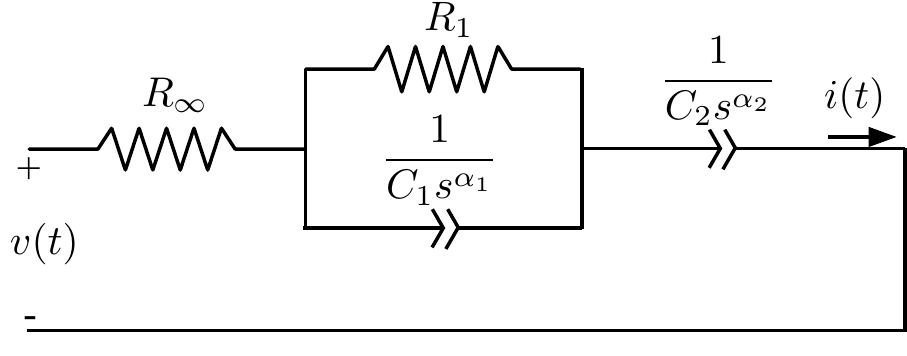}
\caption{The $R_\infty - R_1||\frac{1}{C_1 s^{\alpha_1}} - \frac{1}{C_2 s^{\alpha_2}}$ circuit.}
\label{Fig:fracR-RC-C}
\end{figure}

\begin{figure}
\centering
\includegraphics[scale=0.42]{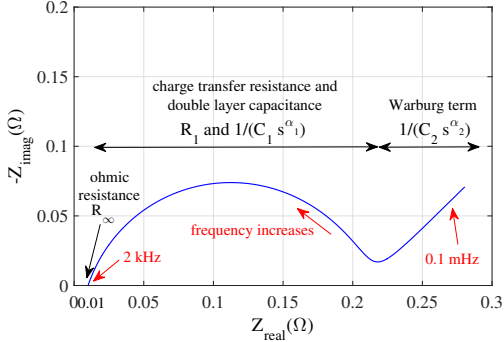}
\caption{Frequency response of the circuit Figure \ref{Fig:fracR-RC-C} for the parameters values given in Eq. \eqref{truevaluesBayes}. } 
\label{Fig:EIS}
\end{figure}


We consider the battery FO model represented in Figure \ref{Fig:fracR-RC-C}, which includes two CPEs. The resistance associated with $C_2$ is open-circuit to model the Warburg term. The model is given by \eqref{SS:d:withnoise} with system matrices \eqref{gen-ecm-ssd} and \eqref{gen-ecm-params} for $n=2$ and $R_2=\infty$. Figure \ref{Fig:EIS} shows the typical frequency response of the circuit which is evident in a large number of electrochemical systems and may be measured using EIS \cite{Barsoukov2005, Orazem2011}. The relation between elements and frequency response have been annotated in the figure. The effect of an additional $R_2$ in parallel with $C_2$ is also investigated later in this study. 
The initial value of the state vector is set to zero: $x_0=[0 ~ 0]$. The parameter vector is
\begin{align*}
\nonumber
\theta &=\left[\begin{array}{cccccc}R_\infty & R_1 & C_1 & C_2 & \alpha_1 &  \alpha_2\end{array}\right].
\end{align*}
We generate synthetic data using the parameter set
\begin{align}
\label{truevaluesBayes}
\theta^\star = \left[\begin{array}{cccccc}0.01 & 0.2 & 3.0 & 400 & 0.8 &  0.5\end{array}\right].
\end{align}

Figure \ref{Fig:EIS} shows the frequency response of the circuit for the above `true' values from 0.1 mHz to 2 kHz. 
\begin{figure*}
\centering
\includegraphics[width=0.48\textwidth]{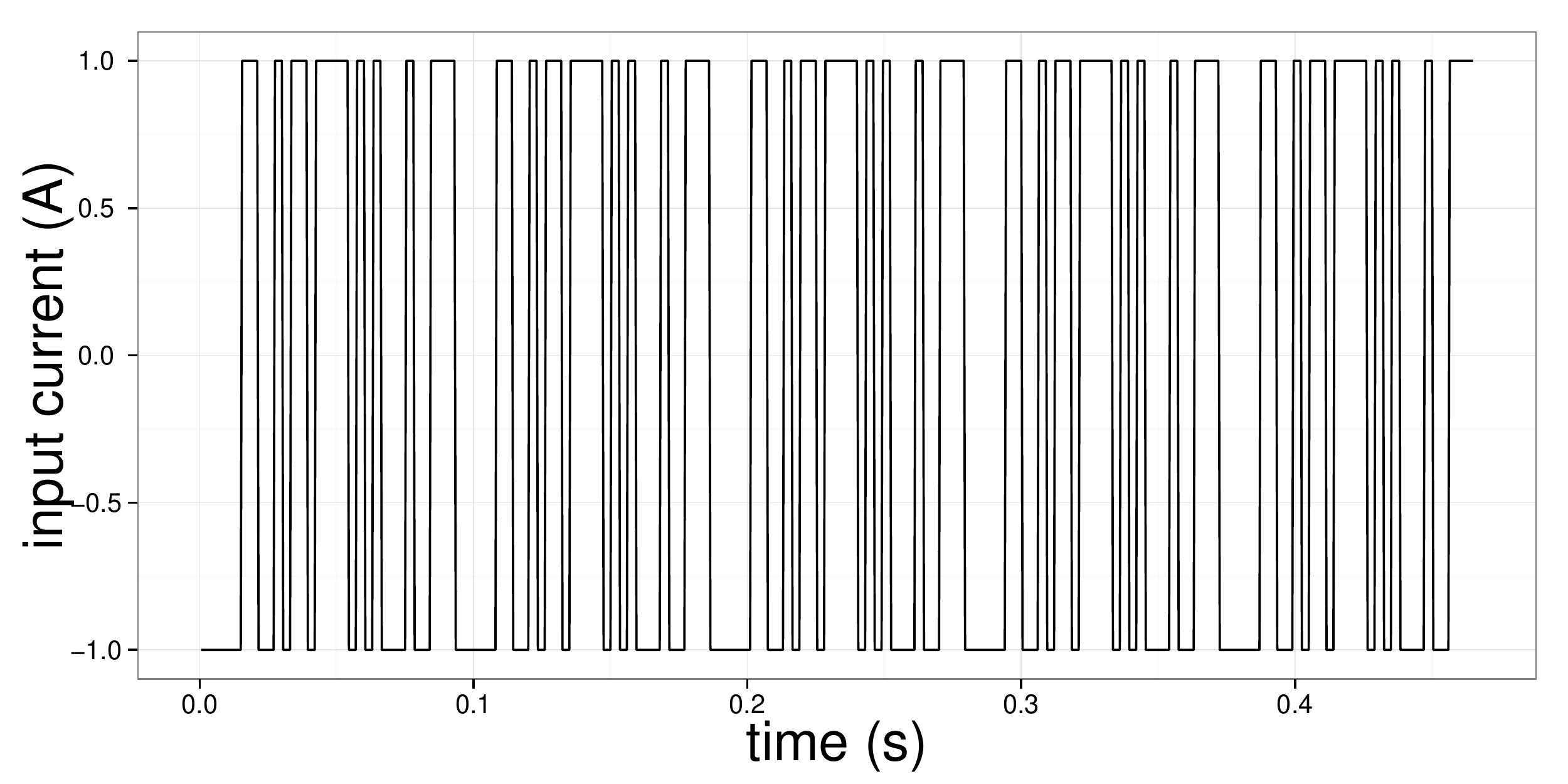}
\includegraphics[width=0.48\textwidth]{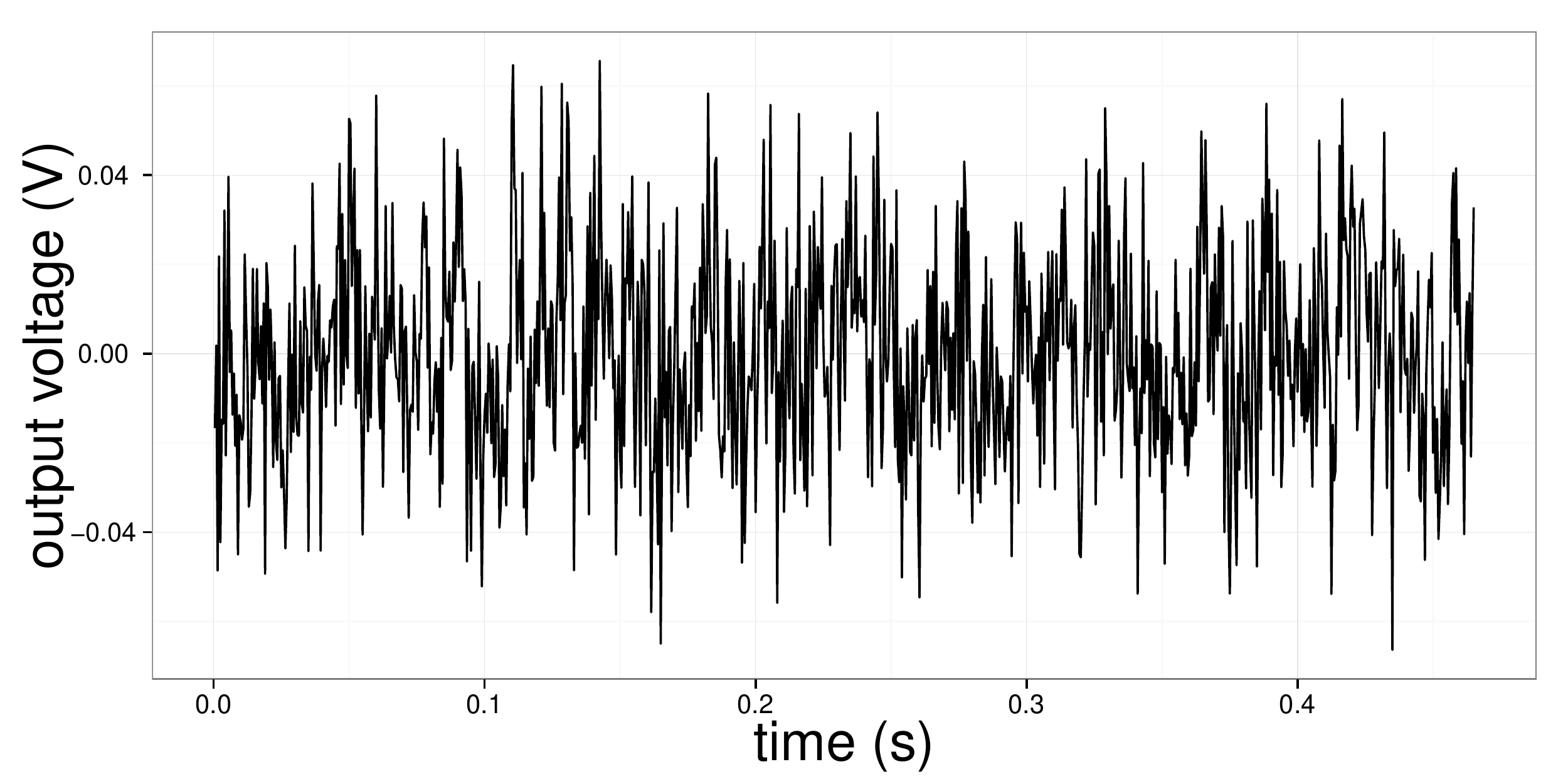}
\caption{Input sequence (left) and observations (right), of length $T=930$, generated from the model of Eq. \eqref{SS:d}. \label{fig:Input-Obs}}
\end{figure*}

The data length is set to $T=930$ samples. The standard deviations are set to $\sigma_{x}=0.002$ and $\sigma_{y}=0.02$. A pseudo-random binary sequence (PRBS) signal between -1 and +1 was generated for $(u_k)_{k\in\mathbb{N}}$, with sampling time $T_s$=0.5 ms. The output voltage $(y_k)_{k\in\mathbb{N}}$ is then generated using the model of Eq. \eqref{SS:d} with the parameter value of Eq. \eqref{truevaluesBayes}. Figure \ref{fig:Input-Obs} shows the input-output data for the base scenario.

%
%

\begin{table}
\caption{\normalsize Range for each parameter.}
\begin{center}
\begin{tabular}{l l  l}
\hline
Parameter& Min & Max \\
\hline

$R_{\infty}~\Omega$ & $0.005$ &   $0.10$  \\
$R_1~~\Omega$         & $0.050$ &  $0.50$ \\
$C_1~~\mbox{Fcm}^{-2}\mbox{s}^{-\alpha_1}$ &  1.00 & 5.00  \\
$C_2~~\mbox{Fcm}^{-2}\mbox{s}^{-\alpha_2}$  & 300  &   500\\
$\alpha_1$, $\alpha_2$   & 0.40  & 1.00 \\
\hline
\end{tabular}
\end{center}
\label{Table:data}
\end{table}

The prior is set to be the uniform distribution over the ranges given in Table \ref{Table:data}. The PMMH is tuned with $N=128$ particles per iteration, and $M=20,000$ iterations. The proposal distribution on $\theta$ is a Gaussian distribution tuned using preliminary runs, in order to match the covariance structure of the posterior distribution (details on this tuning phase are given in Section \ref{sec:algo:pmmh}). We first present the results of the base scenario (Section \ref{sec:results:base}), and then consider various modifications: the number of observations (Section \ref{sec:results:nobs}), the magnitude of the input data (Section \ref{sec:results:mag}), the prior distribution (Section \ref{sec:results:prior}) and the state-output noise ratio of the generated data (Section \ref{sec:results:snr}). For each of the following experiments, five independent runs are performed. The resulting density estimators of each run are overlaid, in order to check the consistency of the method across multiple runs.

\begin{Remark}
  \label{pe-signal-remark}
The ability to identify model parameters depends on the quantity and quality (uncertainties) of the available data. Although there is a rich literature on the experiment design\cite{Ljung1987, Soderstrom1989, Rothenberger2014}, it remains unclear how to generate informative data for parameter estimation of FO systems. In \cite{Alavi2015a}, we applied the persistent excitation concept to the battery Randles circuit models, which are given by ordinary differential equations (ODEs). This method is adopted here and three different persistently exciting PRBS inputs are applied. The authors in \cite{Raue2009} propose a method for the design of periodic excitation signals to identify the battery models given by ODEs. Further investigations in this regard are left open for future research. \hfill{$\square$}
\end{Remark}


\begin{figure*}
\begin{centering}
\includegraphics[width=1\textwidth]{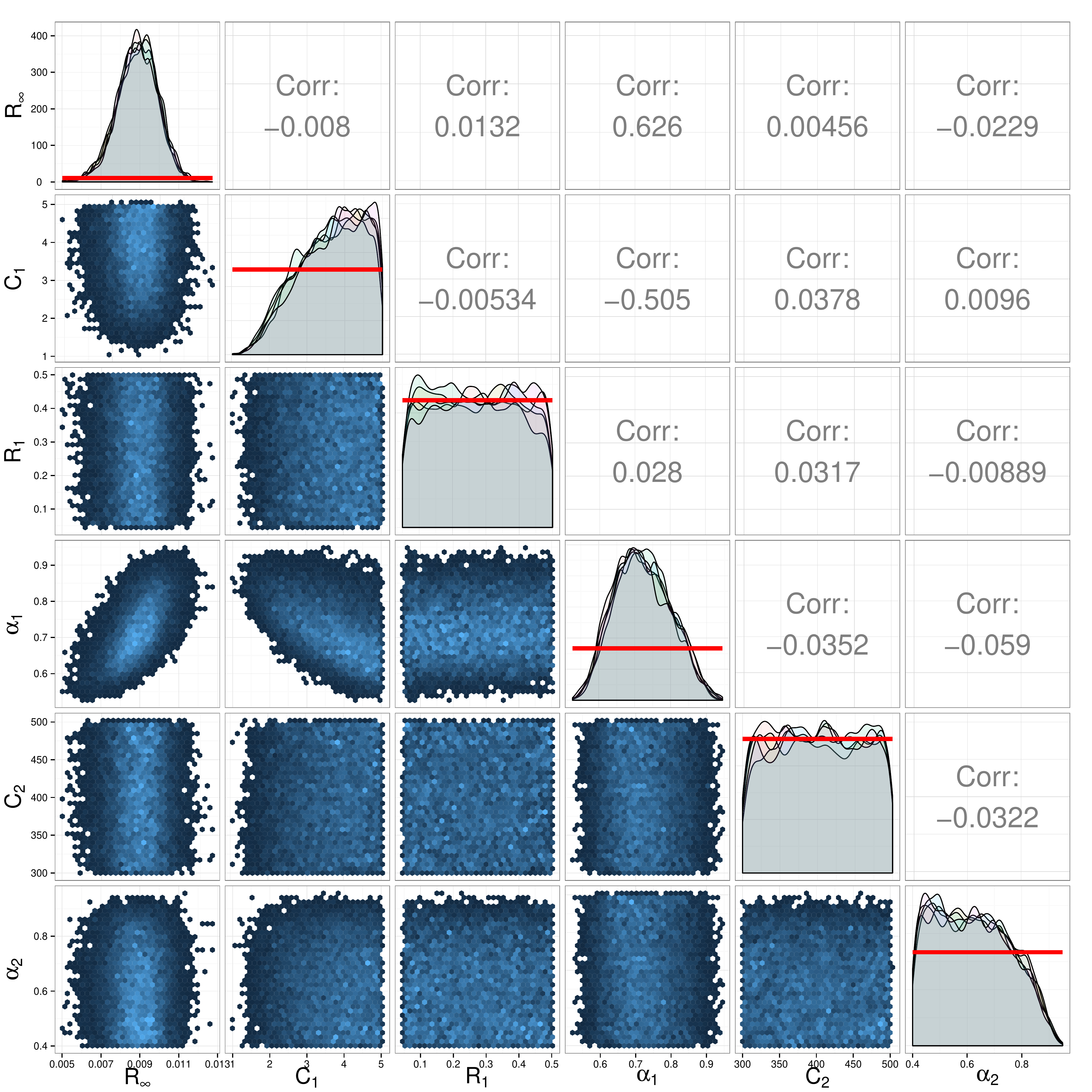}
\par\end{centering}
\protect\caption{Results of the PMMH algorithm for the base scenario, sampling approximately from the posterior distribution for the model \eqref{SS:d}. On the diagonal, the posterior density function of each run is overlaid with the prior distribution (red horizontal lines), which is uniform. On the lower triangle, the five runs are pooled together in a histogram with hexagonal binning. On the upper triangle, the correlation coefficients between pairs of parameters are displayed. \label{fig:pmmhresults:ggpair}}
\end{figure*}

\subsection{Posterior samples in the base scenario \label{sec:results:base}}
The results of the base scenario, computed on five independent runs, are shown in Figure \ref{fig:pmmhresults:ggpair}. The graph shows the approximation of each marginal posterior distribution
using kernel density estimators, computed on each run (diagonals), as well as pairwise histograms with hexagonal binning (lower triangle) of all the runs combined. The pairwise correlations are indicated in the upper triangle.

On the diagonal, we see that the five independent runs are consistent, indicating that the PMCMC method approximates the posterior distribution in a satisfactory way. We can then comment on the posterior itself. We see that for some parameters, such as $R_{\infty}$, the posterior is significantly different from
the prior, represented by horizontal red lines. This indicates that the data are informative on these parameters. Therefore, we can expect that the marginal posterior distributions of these parameters would concentrate around the corresponding components of the data-generating parameter $\theta^\star$ of Eq. \eqref{truevaluesBayes},
when the number of observations goes to infinity. On the other hand, nothing seems to be learned on some other parameters, such as $C_{2}$, for which the posterior resembles the prior. 

Recall that the prior distribution is set to be uniform over the intervals defined by Table \ref{Table:data}. As a result, the posterior distribution is simply proportional to the likelihood of Eq. \eqref{eq:bayesformula}, and thus the mode of the posterior would be precisely the maximum likelihood estimate, under the constraints of Table \ref{Table:data}. However, since the posterior distribution is flat on some parameters (such as $C_{2}$), the mode could be anywhere in the admissible interval for these parameters. Thus, a numerical procedure giving only the maximum likelihood estimate would return any value in that range. 
In \cite{Alavi2015b}, we have seen the similar results that the estimation error of Warburg term using the battery Randles model with standard capacitors, is large. Therefore, the identification of the Warburg term in the battery FO model is similar to the identification of the integral term. 

The Bayesian approach is closer in spirit to integrated likelihood approaches, which are alternatives to profile likelihood. Another advantage of sampling from the posterior distribution is the possibility to investigate correlations between parameters, which are particularly large between $R_{\infty}$ and $\alpha_{1}$, and between $C_{1}$ and $\alpha_{1}$ according to Figure \ref{fig:pmmhresults:ggpair}.

\subsection{Effect of the number of observations \label{sec:results:nobs}}

We now proceed to investigate whether variations in the experiments change the identification of the parameters. We first consider the number of observations. We consider three sets of data of sizes $T=635$, $T=930$ and $T=1890$, generated from three input sequences of these lengths. Instead of showing of all the marginal distributions as in Figure \ref{fig:pmmhresults:ggpair}, we focus on two parameters, $R_{\infty}$ and $C_{2}$, which are respectively easy and challenging to identify (according to Figure \ref{fig:pmmhresults:ggpair}).

\begin{figure*}
\begin{centering}
\includegraphics[width=0.45\textwidth]{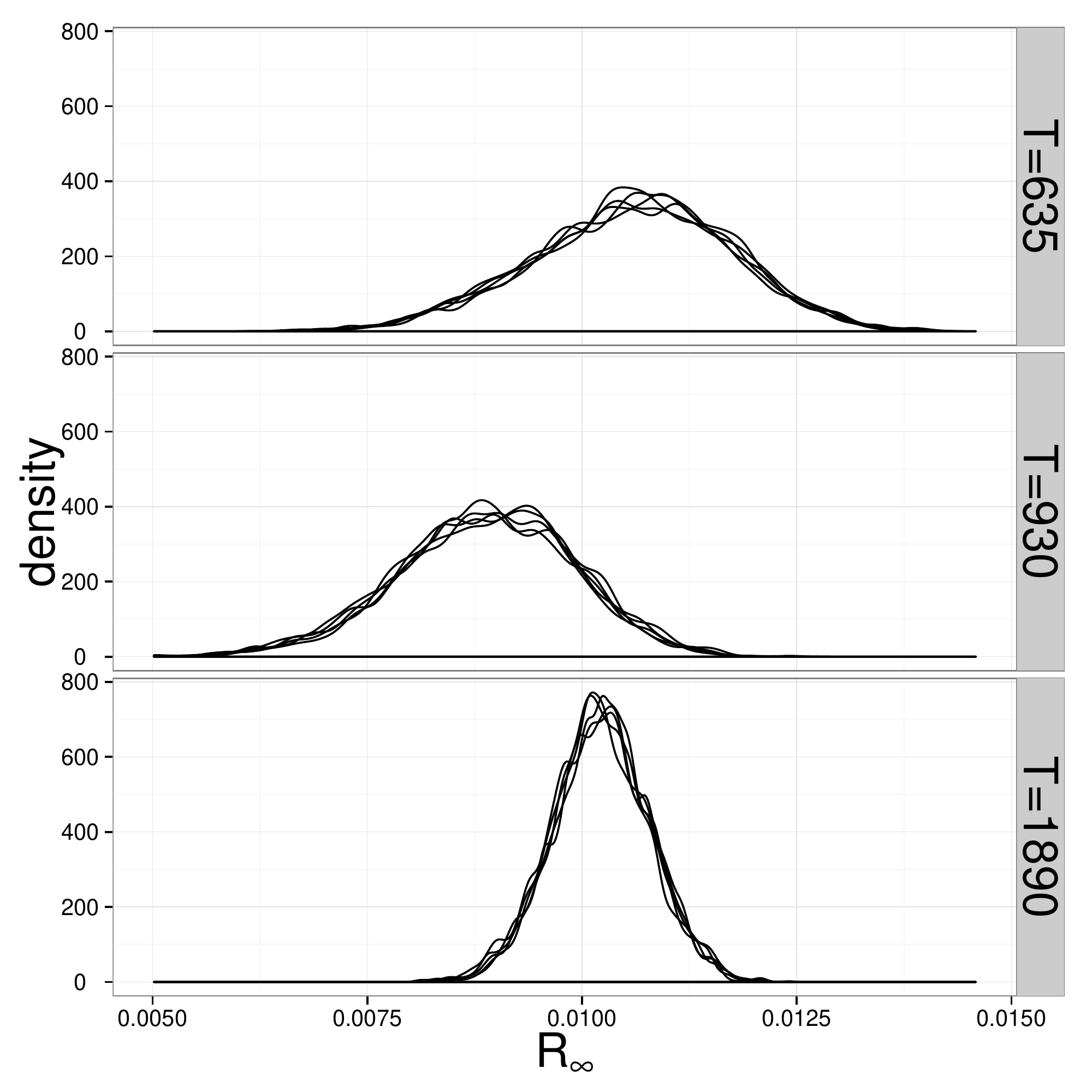}
\includegraphics[width=0.45\textwidth]{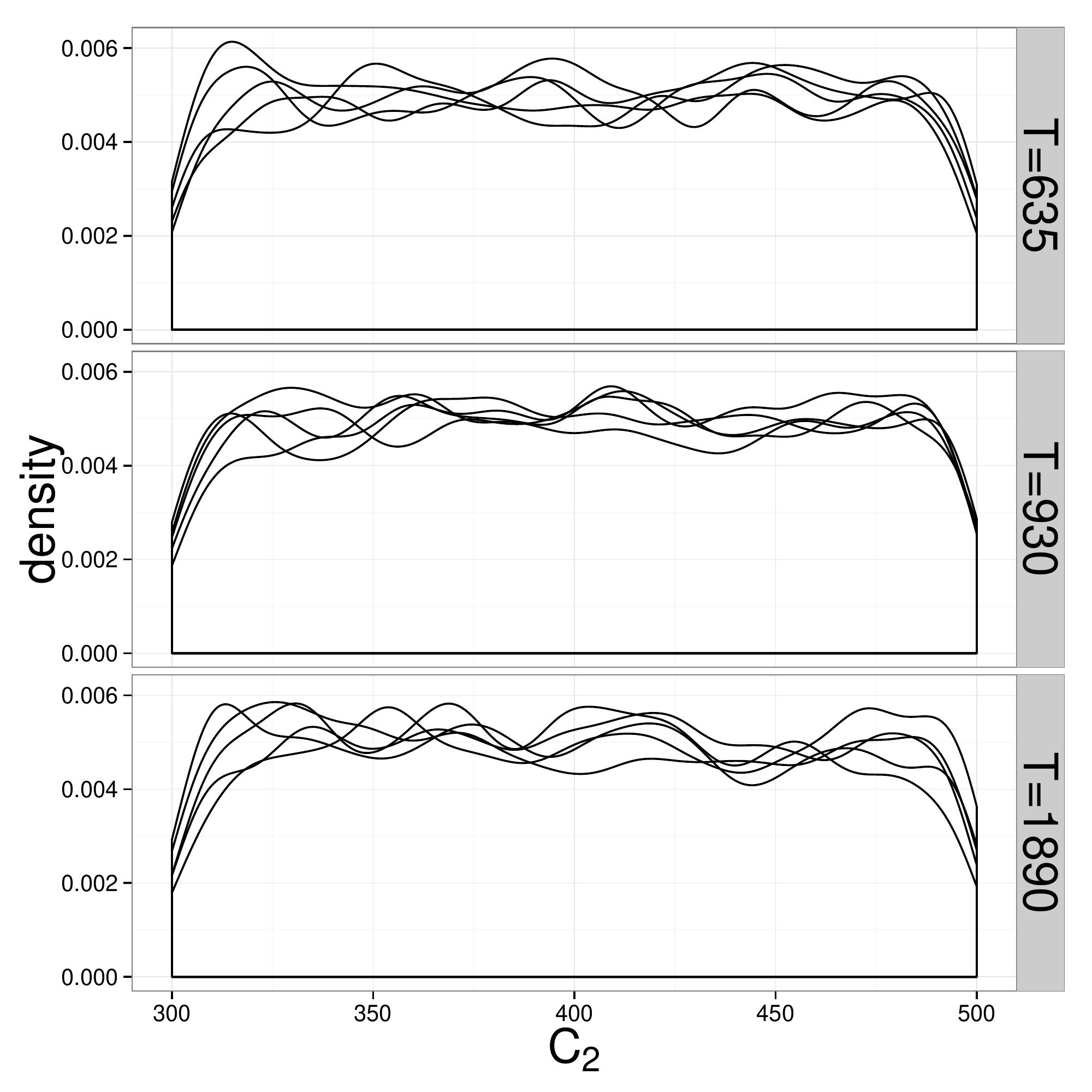}
\par\end{centering}
\protect\caption{Results of the PMMH algorithm, comparing the posterior distributions obtained for $T=635$, $T=930$ and $T=1890$ data points. For some parameters (such as $R_{\infty}$ on the left), adding more data helps the concentration of the posterior distribution, while
for others ($C_{2}$ on the right), increasing data does not seem to have any effect. \label{fig:pmmhresults:datalength}}
\end{figure*}

The results are shown in Figure \ref{fig:pmmhresults:datalength}. We see that adding more data makes the posterior distribution more concentrated
for $R_{\infty}$, and closer to the true value of Eq. \eqref{truevaluesBayes}, whereas it does not seem to have an effect on the distribution of $C_{2}$.

\subsection{Effect of the input magnitude \label{sec:results:mag}}

\begin{figure*}
\begin{centering}
\includegraphics[width=0.45\textwidth]{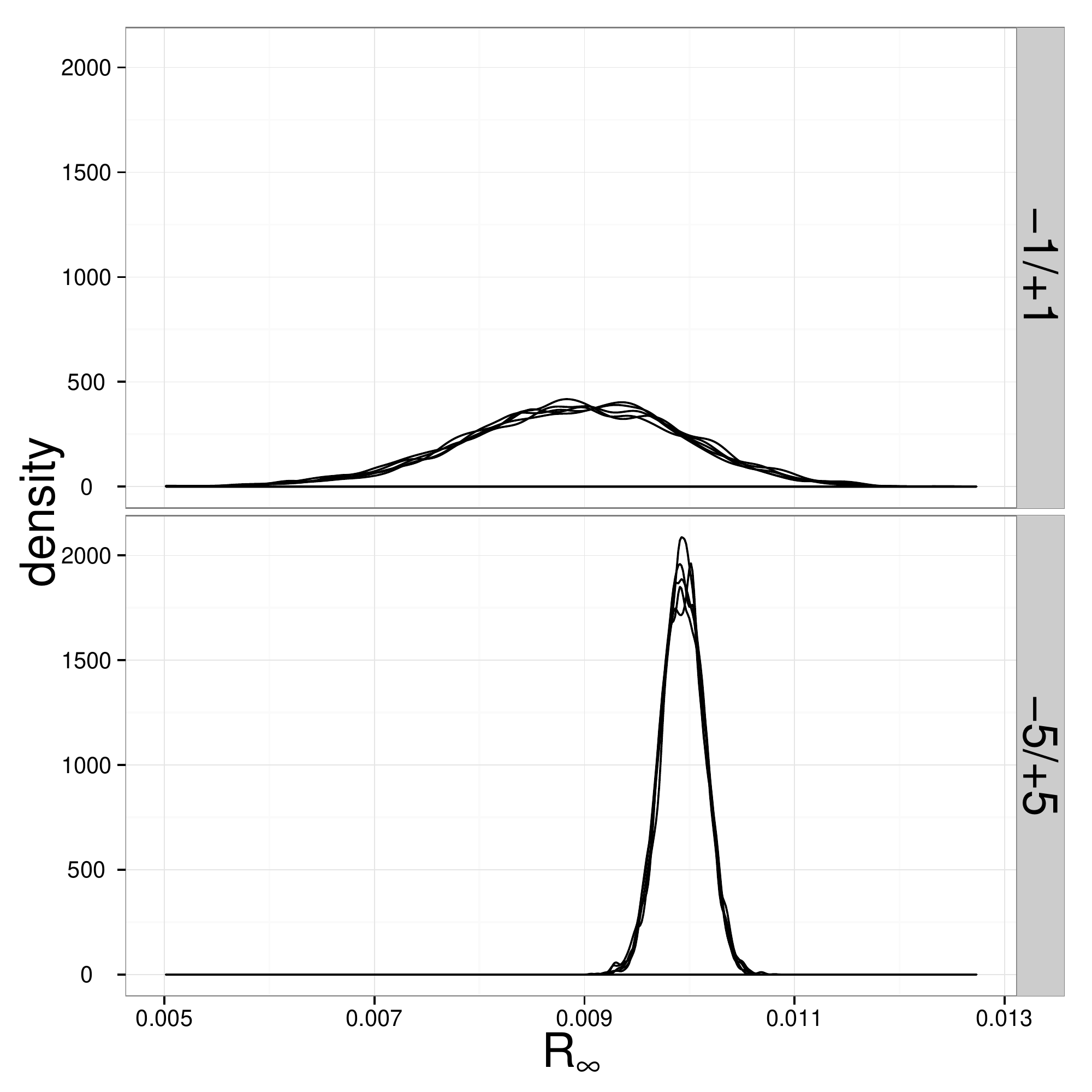}
\includegraphics[width=0.45\textwidth]{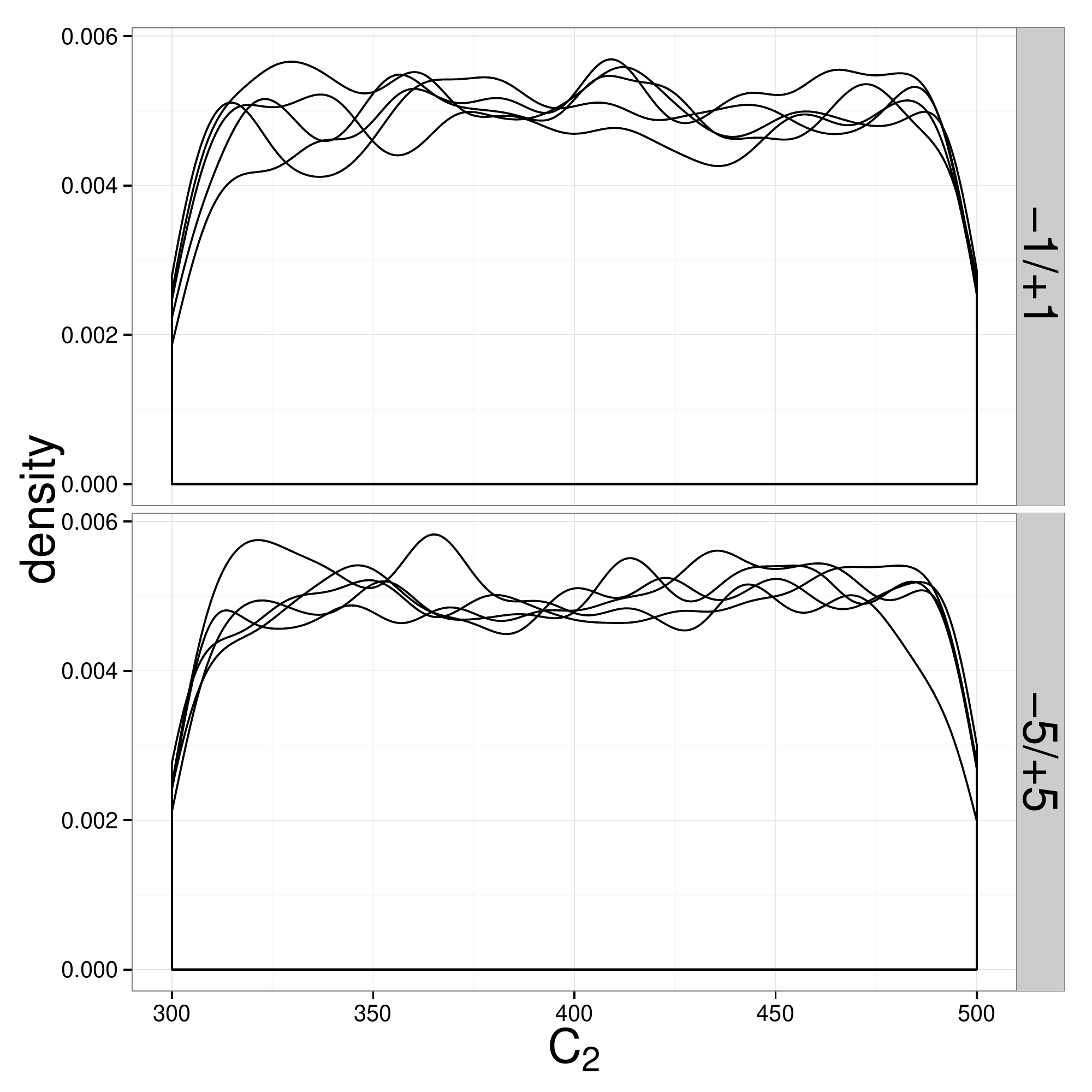}
\par\end{centering}
\protect\caption{Results of the PMMH algorithm, comparing the posterior distributions obtained for inputs of magnitude $1$ (top panels) and of magnitude $5$ (bottom panels). For some parameters (such as $R_{\infty}$ on the left), an input with higher magnitude helps the concentration of the posterior distribution, while for others ($C_{2}$ on the right), changing the magnitude does not seem to have any effect.
\label{fig:pmmhresults:magnitude}}
\end{figure*}

We study the effect of the magnitude of the input. On top of the base scenario, with input data of magnitude $1$, we consider an input sequence of magnitude $5$ (thus, oscillating between $-5$ and $+5$). The results are shown in Figure \ref{fig:pmmhresults:magnitude}. Increasing the magnitude of the input helps identifying $R_{\infty}$, but still does not seem to impact the posterior distribution of $C_{2}$.

\subsection{Effect of the prior distribution \label{sec:results:prior}}

\begin{figure*}
\begin{centering}
\includegraphics[width=0.45\textwidth]{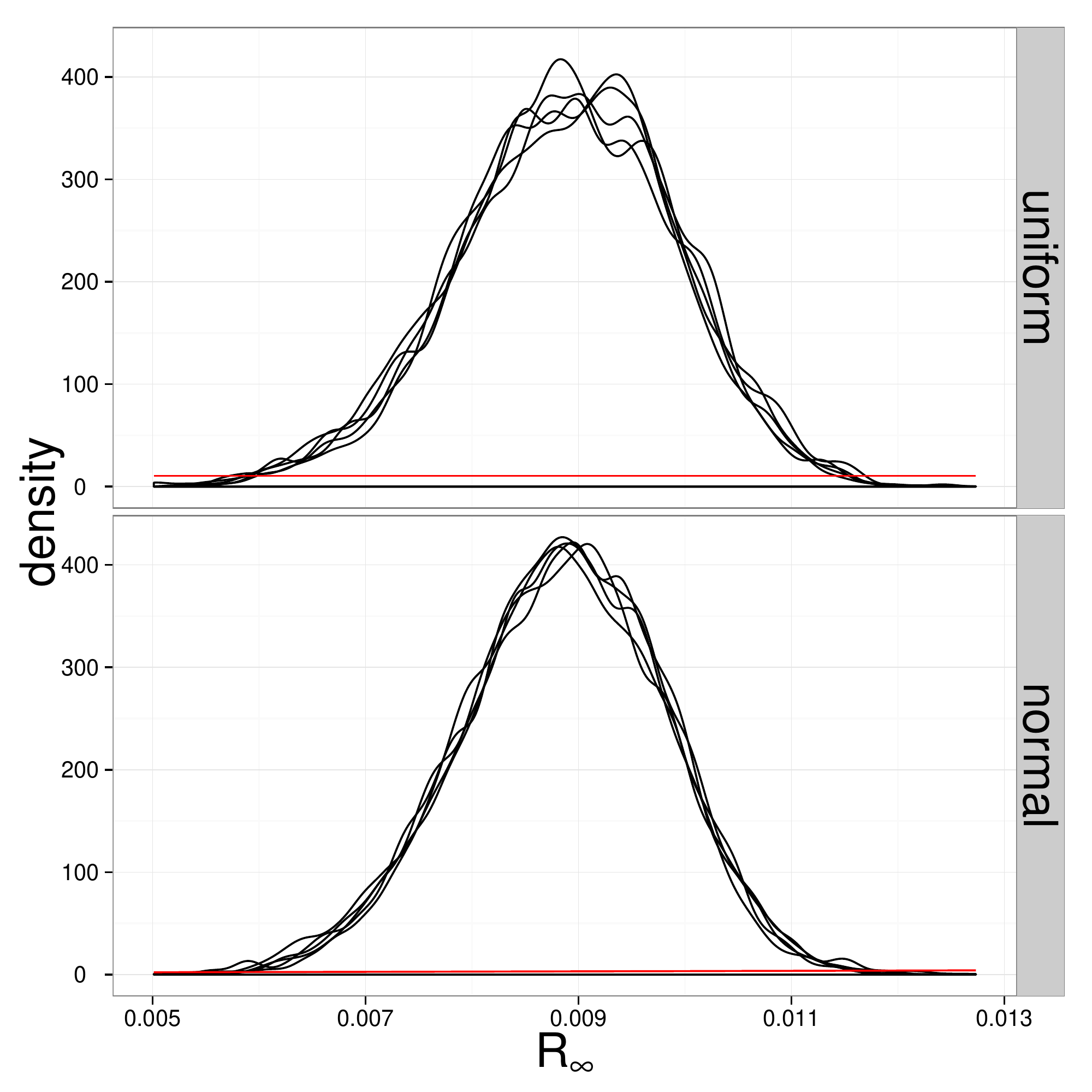}
\includegraphics[width=0.45\textwidth]{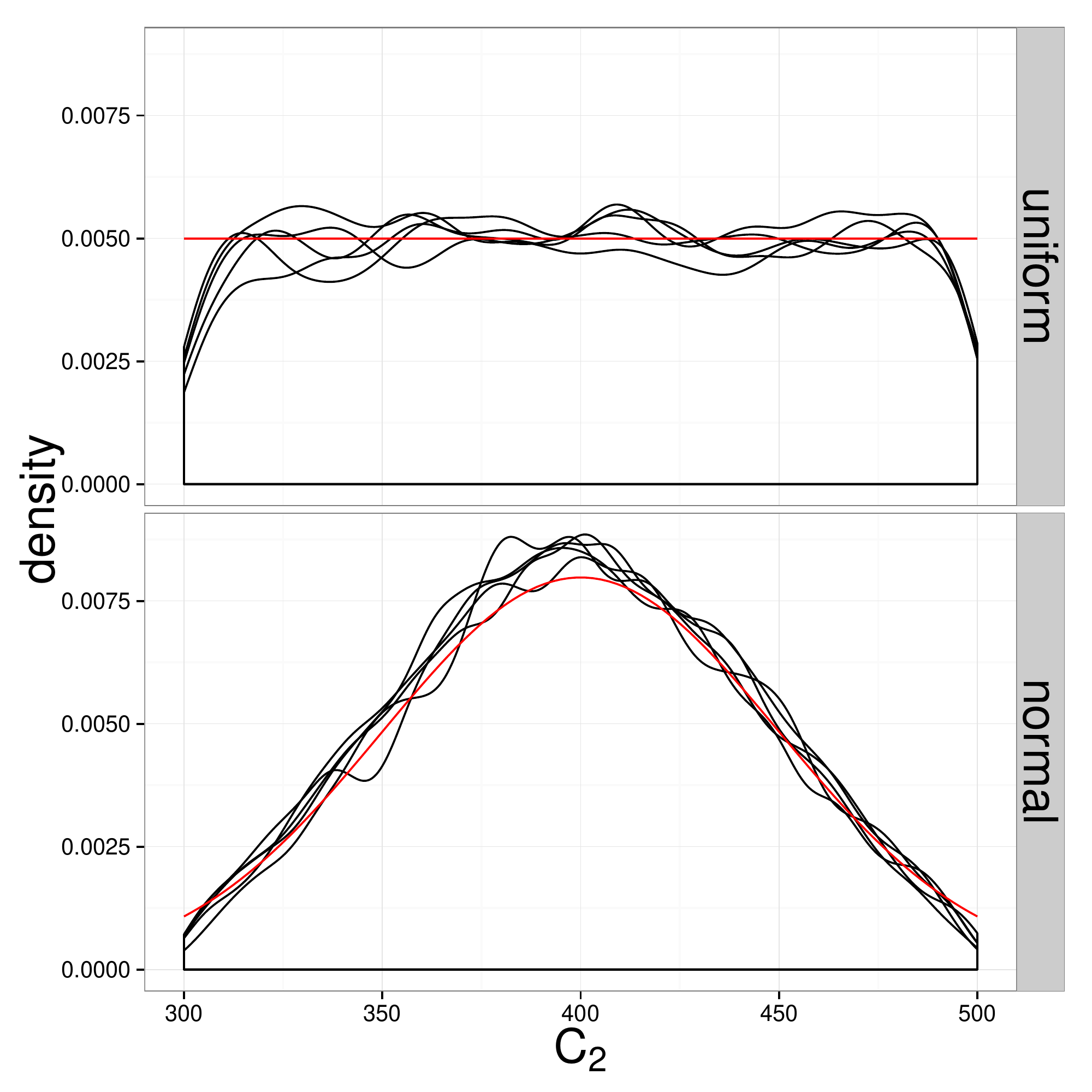}
\par\end{centering}

\protect\caption{Results of the PMMH algorithm, comparing the posterior distributions obtained for two prior distributions: a uniform prior (top panels)
and a normal prior (bottom panels). For some parameters (such as $R_{\infty}$ on the left), the prior distribution has no noticeable effect on the posterior distribution, which is very concentrated with respect to both priors. On some others ($C_{2}$ on the right), the posterior looks very much like the prior and thus the choice of prior has a huge effect on the results. \label{fig:pmmhresults:prior}}
\end{figure*}

We study the effect of the prior distribution. We consider a Gaussian prior, centered in the middle of the range indicated in Table \ref{Table:data}, and with standard deviation one fourth of the range, so that roughly $95\%$ of the Gaussian mass lands in the range; we truncate the Gaussian distribution outside the range. Since some of the parameters are barely identifiable, we expect the choice of prior to have an impact. Figure \ref{fig:pmmhresults:prior} shows that indeed, for some parameters such as $C_{2}$, the posterior is essentially equal to the prior, and hence the choice of prior has a huge impact. For other parameters such as $R_{\infty}$, the choice of prior distribution has no noticeable effect, indicating that the posterior distribution is highly driven by the observations.

\subsection{Effect of the state-output noise ratio \label{sec:results:snr}}

\begin{figure*}
\begin{centering}
\includegraphics[width=0.45\textwidth]{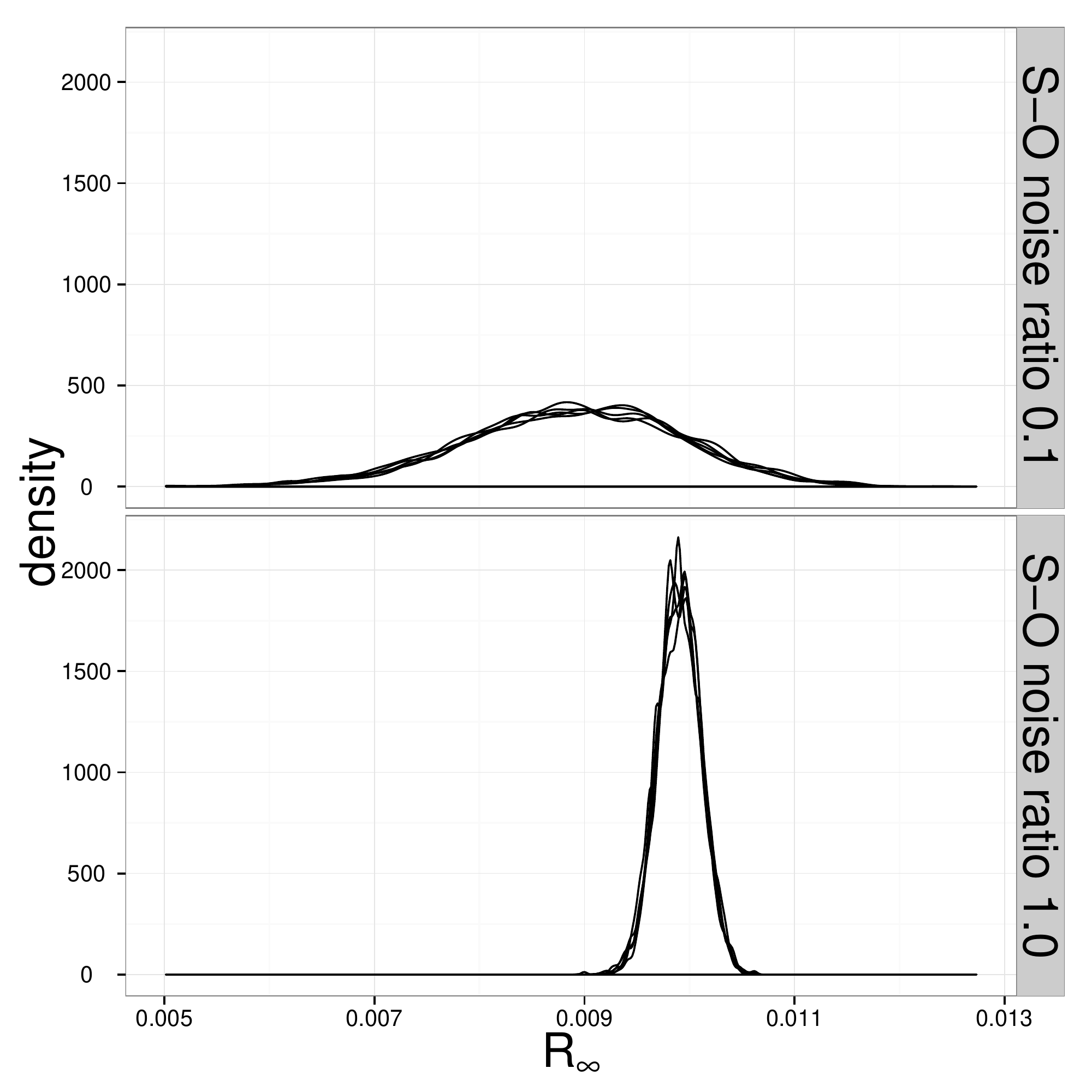}
\includegraphics[width=0.45\textwidth]{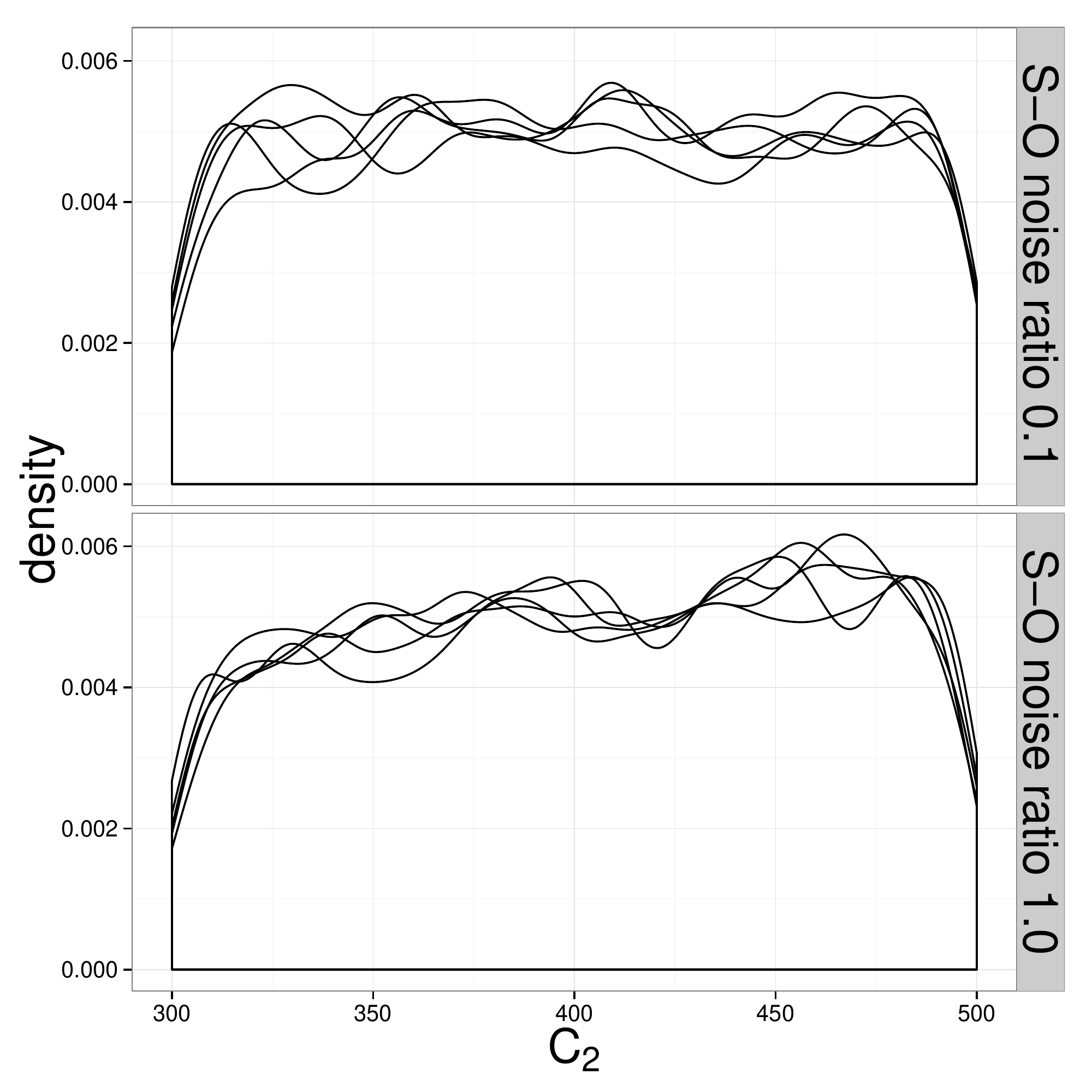}
\par\end{centering}

\protect\caption{Results of the PMMH algorithm, comparing the posterior distributions obtained for state-output noise ratios: a state-output noise ratio of $0.1$ (top panels) and a state-output noise ratio of $1.0$ (bottom panels). For some parameters (such as $R_{\infty}$ on the left), the state-output noise ratio has a huge effect on the posterior, while on some others ($C_{2}$ on the right), the effect is less noticeable, if present at all. 
\label{fig:pmmhresults:snr}}
\end{figure*}

We study the effect of the state to output noise ratio. More precisely, we generate data using $\sigma_{y}=0.002$, instead of using $\sigma_{y}=0.02$ as in the base scenario. We refer to these state-output noise ratios (i.e.\ signal to noise ratios) as $1.0$ and $0.1$ respectively. The results are shown in Figure \ref{fig:pmmhresults:snr}. We see that increasing the state-output noise ratio has more of an effect on some parameters than others. Note that, this time, the posterior distribution of $C_{2}$ seems to deviate from the prior distribution, with a slight inclination of the posterior toward the right-hand side of the interval $\left[300,500\right]$.

In our recent work \cite{Alavi2015b} we showed that non-commensurate FO models are structurally (i.e. assuming noise-free, perfect data) locally identifiable. For the Ohmic resistor $\R_{\infty}$ we showed that it can be globally structurally identifiable. In the present we study the practical identifiability (i.e. assuming real, noisy data) using the framework of Bayesian inference. In particular we note that the term $R_{\infty}$ is also practically globally identifiable (see Figures above).


We believe that identification of the FO Warburg term is similar to the identification of integral terms in ODE systems. 
Closed-loop identification \cite{Nayfeh2004} is a common method for the identification of the integral term. Further research is required to investigate whether the closed-loop identification approach could be applied to the battery non-commensurate FO models.

In the context of Bayesian inference, we suggest other particle MCMC methods such as Particle Gibbs \cite{Andrieu2010,Lindsten2012} or SMC$^2$ \cite{chopin2013smc2} could also be envisioned, with potential savings in computation time.

\section{Conclusions }
In this paper, we have developped a computational approach to Bayesian inference in non-Markovian state-space models, and applied it to the identification of fractional order (FO) systems. The posterior distribution was approximated using a particle marginal Metropolis-Hastings (PMMH) method, and in particular an implementation of particle filters based on a tree representation of the trajectories has been proposed to improve the efficiency of numerical computations. We applied the method to the battery non-commensurate FO model and studied identifiability of the model parameters subject to the choice of prior, excitation signal, data length, and noise ratio. 

In our recent work \cite{Alavi2015b} we showed that non-commensurate FO models are structurally (i.e. assuming noise-free, perfect data) locally identifiable. For the Ohmic resistor $\R_{\infty}$ we showed that it can be globally structurally identifiable. In the present we study the practical identifiability (i.e. assuming real, noisy data) using the framework of Bayesian inference. In particular we note that the term $R_{\infty}$ is also practically globally identifiable (see simulation in Section \ref{sec:batpmcmc}). We believe that identification of the FO Warburg term is similar to the identification of integral terms in ODE systems. 
Closed-loop identification \cite{Nayfeh2004} is a common method for the identification of the integral term. Further research is required to investigate whether the closed-loop identification approach could be applied to the battery non-commensurate FO models.

In the context of Bayesian inference, we suggest other particle MCMC methods such as Particle Gibbs \cite{Andrieu2010,Lindsten2012} or SMC$^2$ \cite{chopin2013smc2} could also be envisioned, with potential savings in computation time.


\section{Appendices}
\subsection{Particle filter \label{sec:algo:particlefilter}}

We describe a general particle filter for the model
of Section \ref{sub:statespacemodels}, in Algorithm \ref{alg:pf}. It relies on a resampling distribution,
denoted by $r(\cdot\mid w_{k-1}^{1:N})$, to sample the ancestor indices
of each generation of particles. There are multiple choices for the
resampling distribution; Section \ref{sub:particlefilters} has described
the multinomial resampling scheme for simplicity. In the numerical results, we use
a more sophisticated method called systematic resampling, described
in Algorithm \ref{alg:systematicresampling}, and explained in \cite{Carpenter1999,Murray2013a}.

\begin{algorithm}
\begin{algorithmic}[1]
	\STATE Sample for each $i\in\left\{ 1,\ldots,N\right\}$: $x_{0}^{i}\sim \mu(\cdot\mid\theta)$, and define $\bar{x}_0^i = x_0^i$.
	\STATE Compute for each $i\in\left\{ 1,\ldots,N\right\}$:
\[w_{0}^{i}=g_{0}\left(y_0\mid x_{0}^{i},\theta\right).\]
    \STATE Compute the likelihood estimator:
\[ \hat{p}(y_0\mid \theta) = \frac{1}{N}\sum_{i=1}^N w_0^i.\]
	\FOR {$k\geq 1$}
		\STATE Draw the ancestors $a_{k-1}^{1:N}$ from the resampling distribution $r(\cdot \mid w_{k-1}^{1:N})$.
		\STATE Sample for each $i\in\left\{1,\ldots,N\right\}$:
				\[x_{k}^{i}\sim q_{k}(\cdot\mid \bar{x}_{0:k-1}^{a_{k-1}^i},\theta).\]
				and define $\bar{x}_{0:k}^i = (\bar{x}_{0:k-1}^{a_{k-1}^i}, x_k^i)$.
		\STATE Compute for each $i\in\left\{1,\ldots,N\right\}$:
\[w_{k}^{i}=\frac{f_{k}(x_{k}^{i}\mid \bar{x}_{0:k-1}^{a_{k-1}^i},\theta)
g_{k}(y_k\mid x_{k}^{i},\theta)}{q_{k}(x_{k}^{i}\mid \bar{x}_{0:k-1}^{a_{k-1}^i},\theta)}.\]
		\STATE Compute the likelihood estimator:
\[ \hat{p}(y_{0:k}\mid \theta) = \hat{p}(y_{0:k-1}\mid \theta) \times \frac{1}{N}\sum_{i=1}^N w_k^i.\]
	\ENDFOR
	\RETURN $\hat{p}(y_{0:k}\mid\theta)$ and $(\bar{x}_{0:k}^i, w_k^i)_{i=1}^N$ for all $k\geq 0$.
\end{algorithmic}

\protect\caption{Particle filter, producing approximate samples $(\bar{x}_{0:k}^{i},w_{k}^{i})$
of the filtering distributions $p(x_{0:k}\mid y_{0:k},\theta)$, and
an estimator $\hat{p}(y_{0:k}\mid\theta)$ of the likelihood $p(y_{0:k}\mid\theta)$,
for all $k\in\mathbb{N}$.\label{alg:pf}}
\end{algorithm}

\begin{algorithm}
\begin{algorithmic}[1]
	\STATE Sample a uniform random variable: $u\sim \mathcal{U}([0,1])$.
	\STATE Normalise the weights: for each $i\in\{1,\ldots,N\}$, $\omega^i = w^i / \sum_{j=1}^N w^j$.
	\STATE Set $\bar{u} = u/N$, $j = 1$, $S_\omega = \omega^1$.
	\FOR {$k=1$ \TO $k= N$}
		\WHILE{$S_\omega < \bar{u}$}
			\STATE Set $j \leftarrow j + 1$.
			\STATE Set $S_\omega \leftarrow S_\omega + \omega^j$.
		\ENDWHILE
		\STATE Set $a^k = j$.
		\STATE Set $\bar{u} \leftarrow \bar{u} + 1/N$.
	\ENDFOR
	\RETURN $a^{1:N}$.
\end{algorithmic}

\protect\caption{Systematic resampling, returning a vector of ancestors $a^{1:N}$
given a vector of weights $w^{1:N}$.\label{alg:systematicresampling}}
\end{algorithm}

\subsection{Locally optimal proposal for the battery model \label{sec:localproposal}}

In the models considered in the article, the locally optimal proposal
can be implemented. Indeed the state $x_{k}\in\mathbb{R}^{2}$ follows
a Gaussian distribution centered at a function of the past: $x_{k}\sim\mathcal{N}(\varphi_{k},\Sigma_{x})$,
where $\varphi_{k}=\left(\varphi_{k,1},\varphi_{k,2}\right)^{T}$
is computed as a deterministic function of the past trajectory $x_{0:k-1}$ and of the parameter,
and $\Sigma_{x}=\text{diag}(\sigma_{x}^{2},\sigma_{x}^{2})$. Furthermore
we have $y_{k}\sim\mathcal{N}(c_{k}+x_{k,1}+x_{k,2},\sigma_{y}^{2})$,
where $c_{k}$ is a real value that can be computed given the parameter 
value. As a result, we can compute
$y_{k}\mid x_{0:k-1}\sim\mathcal{N}(\zeta_{k},2\sigma_{x}^{2}+\sigma_{y}^{2})$,
where $\zeta_{k}=c_{k}+\varphi_{k,1}+\varphi_{k,2}$ and $x_{k}\mid x_{0:k-1},y_{k}$
has a Gaussian distribution with
\begin{align*}
\text{mean = } & \begin{pmatrix}\varphi_{k,1}+\frac{\sigma_{x}^{2}}{2\sigma_{x}^{2}+\sigma_{y}^{2}}\left(y_{k}-\zeta_{k}\right)\\
\varphi_{k,2}+\frac{\sigma_{x}^{2}}{2\sigma_{x}^{2}+\sigma_{y}^{2}}\left(y_{k}-\zeta_{k}\right)
\end{pmatrix},\\
\text{and variance = } & \begin{pmatrix}\sigma_{x}^{2}-\frac{\sigma_{x}^{4}}{2\sigma_{x}^{2}+\sigma_{y}^{2}} & -\frac{\sigma_{x}^{4}}{2\sigma_{x}^{2}+\sigma_{y}^{2}}\\
-\frac{\sigma_{x}^{4}}{2\sigma_{x}^{2}+\sigma_{y}^{2}} & \sigma_{x}^{2}-\frac{\sigma_{x}^{4}}{2\sigma_{x}^{2}+\sigma_{y}^{2}}
\end{pmatrix}.
\end{align*}

\subsection{Particle Marginal Metropolis--Hastings \label{sec:algo:pmmh}}

The Particle Marginal Metropolis--Hastings algorithm of \cite{Andrieu2010}
is described in Algorithm \ref{alg:pmmh}. It requires a few choices
of tuning parameters from the user, which we now detail.

First, it requires a particle filter yielding likelihood estimates
$\hat{p}(y_{0:T}\mid\theta)$ for all $\theta$, as in Section \ref{sec:algo:particlefilter}. We have described
a choice of proposal distribution in Section \ref{sec:localproposal}, but one still needs
to choose a number of particles $N$. We adopt the following strategies,
inspired by the ``conditional acceptance rate'' of \cite{Murray2013}.
For a sample $\theta$ from the prior distribution, we run the particle
filter $n$ times, given that same $\theta$. This yields a sequence
of estimates $Z_{1},\ldots,Z_{n}$ of $p(y_{0:T}\mid\theta)$. We
then mimic a Metropolis--Hastings scheme with these estimates, as
follows:
\begin{itemize}
\item Start with $Z_{\text{current}}=Z_{1}$.
\item For all $j=2,\ldots,n$,

\begin{itemize}
\item compute the probability
\[
\alpha=\min\left(1,\frac{Z^{j}}{Z_{\text{current}}}\right),
\]

\item with probability $\alpha$, accept $Z_{\text{current}}=Z^{j}$, otherwise
keep $Z_{\text{current}}$ unchanged.
\end{itemize}
\item Report the average number of times that a new $Z$ was accepted.
\end{itemize}
The average number of acceptances would be $100\%$ if the estimates
were all perfect evaluations of $p(y_{0:T}\mid\theta)$. On the other
hand, the number of acceptances is very low if the variance of the
estimates is very large. Thus, starting from a small number of particles
$N$, we increase the number $N$ until the average number of acceptances
is considered to be enough, for instance $10\%$. We can perform this
procedure on a few parameters $\theta$ in parallel, to check the
variability of the conditional acceptance rate across the parameter
space.

Once the particle filter is tuned, we need to choose a proposal distribution
to draw the candidate parameters: $\theta^{\star}\sim q_{\theta}(\cdot|\theta)$.
A standard choice consists in using a Gaussian random walk:
\[
\theta^{\star}\sim\mathcal{N}\left(\theta,\Sigma\right),
\]
for which we still need to choose the covariance matrix. One approach
consists in changing $\Sigma$ along the run of the algorithm
\cite{Andrieu2008, Roberts2009,Peters2010}. We take
a simpler approach, using a preliminary run.
\begin{itemize}
\item A first PMMH run of $5,000$ iterations is performed, using a diagonal
matrix $\Sigma_{\text{prior}}$. The diagonal elements are chosen
to match the standard deviation of the prior distribution. In our
numerical experiments, the obtained acceptance rate is then low but
non-zero.
\item We discard the first $2,500$ iterations, and compute the covariance
matrix $\hat{\Sigma}$ of the last $2,500$ samples, which is a crude
estimate of the covariance of the posterior distribution. We also
retain the last sample of these $2,500$ samples.
\item The PMMH algorithm is then run again, for $20,000$, using $\Sigma=\hat{\Sigma}$
and starting the chain from the last sample of the preliminary run.
We observe that the acceptance rate during this final stage is higher
than during the first stage, which indicate an improvement in using
$\hat{\Sigma}$ instead of $\Sigma_{\text{prior}}$.
\end{itemize}
\begin{algorithm}
\begin{algorithmic}[1]
	\STATE Start from a draw from the prior distribution: $\theta^{0}\sim p(\cdot)$.
	\STATE Run a particle filter with $N$ particles to obtain $\hat{p}(y_{0:T}\mid \theta^{0})$.
	\FOR {$t=1$ \TO $t=M$}
		\STATE Propose $\theta^\star \sim q_\theta(\cdot | \theta^{t-1})$.
		\STATE Run a particle filter with $N$ particles to obtain $\hat{p}(y_{0:T}\mid \theta^{\star})$.
		\STATE Compute:
\[\alpha=\text{min}\left(1,
\frac{p(\theta^{\star})\hat{p}(y_{0:T}\mid\theta^{\star})}{p(\theta^{t-1})\hat{p}(y_{0:T}\mid\theta^{t-1})}\frac{q(\theta^{t-1}\mid\theta^{\star})}{q(\theta^{\star}\mid\theta^{t-1})}\right).\]
       \STATE Set $(\theta^{t},\hat{p}(y_{0:T}\mid\theta^t)) = \begin{cases}
                   (\theta^{\star}, \hat{p}(y_{0:T}\mid\theta^{\star})) \quad \text{with probability }\alpha,\\          					(\theta^{t-1}, \hat{p}(y_{0:T}\mid\theta^{t-1})) \quad \text{with probability } 1 - \alpha.\end{cases}$
	\ENDFOR
	\RETURN $(\theta^{t})_{t=1}^{M}$.
\end{algorithmic}

\protect\caption{Particle Marginal Metropolis--Hastings algorithm, producing a sample
from $p(\theta\mid y_{0:T})$.\label{alg:pmmh}}

\end{algorithm}

\section*{Acknowledgements}
Authors would like to acknowledge fundings from the UK's Engineering and Physical Sciences Research Council (EPSRC), S.M.M. Alavi and D.A. Howey under grant EP/K503769/1, A. Mahdi and S.J. Payne under grant EP/K036157/1 and P.E. Jacob under grant EP/K009362/1.

\bibliographystyle{IEEEtran}
\bibliography{Final_Bayesian_frac_sysV3.bbl}
\end{document}